\documentclass[12pt]{article}
\usepackage{a4wide}
\usepackage{amssymb}
\usepackage{amsmath}
\usepackage[russian]{babel}

  \chardef\No=242
\newtheorem{theorem}{Theorem}{}
\newtheorem{definition}{Definition}{}
\newtheorem{corollary}{Corollary}{}
\newtheorem{lemma}{Lemma}{}
\newtheorem{remark}{Remark}{}

\begin{document}

\begin{center}
{\bf On a two-weight criteria for multidimensional
Hardy type operator\\ in $p$-convex Banach function spaces and some application}
\end{center}

\begin{center}
ROVSHAN A.BANDALIEV
\end{center}

{\small  ABSTRACT. The main goal of this paper is to prove a two-weight criteria for multidimensio-nal Hardy type operator from weighted Lebesgue spaces into $p$-convex weighted Banach function spaces. Analogously problem for the dual operator is considered. As application we prove a two-weight criteria for boundedness of multidimensional geometric mean operator and sufficient condition on the weights for boundedness of certain sublinear operator from weighted Lebesgue spaces into weighted Musielak-Orlicz spaces.}

\vspace{2mm}


\vspace{4mm}
\begin{center}
1. {\bf  Introduction.}
\end{center}

The investigation of Hardy operator in weighted Banach function spaces (BFS) have recently history. The goal of this investigations were closely connected with the found of criterion on the geometry and on the weights of BFS for validity of boundedness of Hardy operator in BFS. Characterization of the mapping properties such as boundedness and
compactness were considered in the papers [7], [8], [12], [30] and e.t.c. More precisely, in [7] and [8] were considered the boundedness of certain integral operator in ideal Banach spaces. In [12] was proved the boundedness of Hardy operator in Orlicz spaces. Also, in [30] the compactness and measure of non-compactness of Hardy type operator in Banach function spaces was proved. But in this paper we consider the boundedness of Hardy operator in $p$-convex Banach function spaces and find a new type criterion on the weights for validity of Hardy inequality.
Note that the notion of BFS was introduced in [32]. In particular, the weighted Lebesgue spaces, weighted Lorentz spaces, weighted variable Lebesgue spaces, variable Lebesgue spaces with mixed norm, Musielak-Orlicz spaces and e.t.c. is BFS.

In this paper, we establish an integral-type necessary and sufficient condition on weights, which provides the boundedness of the multidimensional Hardy type operator from weighted Lebesgue spaces into $p$-convex weighted BFS.
We also investigate the corresponding problems for the dual operator. It is well known that the classical two weight
inequality for geometric mean operator is closely connected with the one-dimensional Hardy inequality (see [20]). Analogously,the P\'{o}lya-Knopp type inequalities with multidimensional geometric mean operator are connected with the multidimensional Hardy type operator. Therefore, in this paper, as an application of Hardy inequality we prove the boundedness of multidimensional geometric mean operator and boundedness of certain sublinear operator from weighted Lebesgue spaces into weighted Musielak-Orlicz spaces.

\vspace{3mm}
\begin{center}
2. {\bf Preliminaries}
\end{center}

Let $(\Omega,\, \mu)$ be a complete $\sigma$-finite measure space. By $L_0= L_0(\Omega, \mu)$ we denote the
collection of all real-valued $\mu$-measurable functions on $\Omega.$

\begin{definition}{[32, 29, 6]}
We say that real normed space $X$ is a Banach function space (BFS) if:

(P1)$\;$ the norm $\|f\|_X$ is defined for every $\mu$-measurable function $f,$ and $f\in X$ if and only if $\|f\|_X< \infty;$ $\|f\|_X= 0$ if and only if $f= 0$ a.e. ;

(P2)$\;$  $\|f\|_X= \||f|\|_X$ for all $f\in X;$

(P3)$\;$ if $0\le f\le g$ a.e., then $\|f\|_X\le \|g\|_X;$

(P4)$\;$ if $0\le f_n\uparrow f\le g$ a.e., then $\left\|f_n\right\|_X \uparrow\|f\|_X$ (Fatou property);

(P5)$\;$ if $E$ is a measurable subset of $\Omega$ such that $\mu(E)< \infty,$ then $\left\|\chi_E\right\|_X< \infty,$ where $\chi_E$ is the characteristic function of the set $E;$

(P6)$\;$ for every measurable set $E\subset \Omega$ with $\mu(E)< \infty,$ there is a constant $C_E> 0$ such that  $\int_E f(x)\, dx\le C_E\,\|f\|_X.$
\end{definition}

Given a BFS $X$ we can always consider its associate space $X'$ consisting of those $g\in L_0$ that $f\cdot g\in L_1$ for every $f\in X$ with the usual order and the norm $\|g\|_{X'}= \sup\left\{\|f\cdot g\|_{L_1}:\;\|g\|_{X'}\le 1\right\}.$ Note that $X'$ is a BFS in $(\Omega,\,\mu)$ and a closed norming subspaces.

Let $X$ be a BFS and $\omega$ be a weight, that is, positive Lebesgue measurable and a.e. finite functions on $\Omega.$ Let $X_{\omega}=\left\{f\in L_0:\; f\omega\in X\right\}.$ This space is a weighted BFS equipped with the norm $\|f\|_{X_\omega}= \|f\,\omega\|_{X}.$ (For more detail and proofs of results about BFS we refer the reader to [6] and [29].)

Note that the notion of BFS was introduced in [32].

Let us recall the notion of $p$-convexity and $p$-concavity of BFS's.
\begin{definition}{[43]}
Let $X$ is a BFS. Then $X$ is called $p$-convex for $1\le p\le \infty$ if there exists a constant $M> 0$ such that for all $f_1,\ldots,f_n\in X$
$$
\left\|\left(\sum\limits_{k= 1}^n |f_k|^p\right)^{\frac 1p}\right\|_X\le M\, \left(\sum\limits_{k= 1}^n \left\|f_k\right\|_X^p\right)^{\frac 1p} \;\; \mbox{if}\;\; 1\le p< \infty,
$$
or $\displaystyle{\left\|\sup\limits_{1\le k\le n} \left|f_k\right|\right\|_X\le M\,\max\limits_{1\le k\le n} \left\|f_k\right\|_X}$ if $p= \infty.$ Similarly $X$ is called $p$-concave for $1\le p\le \infty$ if there exists a constant $M> 0$ such that for all $f_1,\ldots,f_n\in X$
$$
\left(\sum\limits_{k= 1}^n \left\|f_k\right\|_X^p\right)^{\frac 1p} \le M\,\left\|\left(\sum\limits_{k= 1}^n |f_k|^p\right)^{\frac 1p}\right\|_X \;\; \mbox{if}\;\; 1\le p< \infty,
$$
or $\displaystyle{\max\limits_{1\le k\le n} \left\|f_k\right\|_X\le M\,\left\|\sup\limits_{1\le k\le n} \left|f_k\right|\right\|_X}$ if $p= \infty.$
\end{definition}

\begin{remark}
Note that the notions of $p$-convexity, respectively $p$-concavity are closely related to the notions of upper $p$-estimate (strong $\ell_p$- composition property), respectively lower $p$-estimate (strong $\ell_p$-decomposition property) as can be found in [29].
\end{remark}

Now we reduce some examples of $p$-convex and respectively $p$-concave BFS.
Let $R^{n}$ be the $n$-dimensional Euclidean space of points $x=\left(x_{1},..., x_{n}\right)$
and let $\Omega$ be a Lebesgue measurable subset in $R^n$ and $\displaystyle{|x|= \left(\sum\limits_{i= 1}^n
x_i^2\right)^{1/2}}.$ The Lebesgue measure of a set $\Omega$ will be denoted by $|\Omega|.$ It is well known that
$\displaystyle{|B(0,1)|= \frac {\pi^{\frac n2}}{\Gamma\left(\frac n2+ 1\right)}},$ where $B(0,1)= \left\{x:\,x\in R^n\right.;$ \linebreak $\left.|x|< 1\right\}.$
\\
{\bf Example 1.1.} Let $1\le q\le \infty$ and $X= L_q.$ Then the space $L_q$ is $p$-convex ($p$-concave) BFS if and only if $1\le p\le q\le \infty$ ($1\le q\le p\le \infty.$)

The proof implies from usual Minkowski inequality in Lebesgue spaces.\\
{\bf Example 1.2.} The following Lemma shows that the variable Lebesgue spaces $L_{q(y)}(\Omega)$ is $p$-convex BFS.
\begin{lemma}{[1]} Let $1\le p\le q(x)\le \overline q< \infty$ for all $y\in \Omega_2\subset R^m.$
Then the inequality
$$
\left\|\|f\|_{L_{p}\left(\Omega_1\right)}\right\|_{L_{q(\cdot)}\left(\Omega_2\right)}\le
C_{p,q}\,\left\|\|f\|_{L_{q(\cdot)}\left(\Omega_2\right)}\right\|_{L_{p}\left(\Omega_1\right)}
$$
is valid, where $\displaystyle{C_{p, q}= \!\!\left(\left\|\chi_{\Delta_1}\right\|_\infty+ \left\|\chi_{\Delta_2}\right\|_\infty+  p\left(\frac 1{\underline q}- \frac{1}{\overline q}\right)\right)\!\!\!\left(\left\|\chi_{\Delta_1}\right\|_\infty+ \left\|\chi_{\Delta_2}\right\|_\infty\right)},$
$\!\underline q= \!\!\mbox {ess}\,\!\inf\limits_{\Omega_2}\! q(x),$ $\!\!\!\!\overline q= \mbox {ess}\,\sup\limits_{\Omega_2} q(x),$ $\Delta_1=\left\{(x, y)\in\Omega_1\times \Omega_2:\, q(y)= p\right\},$ $\Delta_2= \Omega_1\times \Omega_2\setminus \Delta_1$ and $f:\Omega_1\times \Omega_2\rightarrow R$ is any measurable function such that
$$
\left\|\|f\|_{L_p\left(\Omega_1\right)}\right\|_{L_{q(\cdot)}\left(\Omega_2\right)}= \inf \left\{ \mu>
0:\;\;\int\limits_{\Omega_2}\left(\frac{\|f(\cdot,y)\|_{L_p\left(\Omega_1\right)}}
{\mu}\right)^{q(y)}\,dy\le 1\right\}< \infty
$$
and $\displaystyle{\|f(\cdot,y)\|_{L_p\left(\Omega_1\right)}= \left(\int_{\Omega_1} |f(x, y)|^p\,dx\right)^{1/p}}.$
\end{lemma}

Analogously, if $1\le q(x)\le p< \infty,$ then $L_{q(x)}(\Omega)$ is $p$-concave BFS.

\begin{definition}{[38, 15].} Let $\Omega \subset R^n$ be a Lebesgue measurable set. A real function
$\varphi: \Omega\times [0, \infty)\mapsto [0, \infty)$ is called a generalized $\varphi$-function if it satisfies:

a)\, $\varphi(x,\,\cdot)$ is a $\varphi$-function for all $x\in \Omega,$ i.e., $\varphi(x,\,\cdot): [0, \infty)\mapsto [0, \infty)$ is convex and satisfies $\varphi(x,\,0)= 0,$ $\lim\limits_{t\to +0} \varphi(x,\,t)= 0 ;$

b)\, $\psi: x\mapsto \varphi(x,\,t)$ is measurable for all $t\ge 0.$
\end{definition}

If $\varphi$ is a generalized $\varphi$-function on $\Omega,$ we shortly write $\varphi\in \Phi.$
\begin{definition} {[38, 14].} Let $\varphi\in \Phi$  and be $\rho_{\varphi}$  defined by the expression
$$
\rho_{\varphi}(f):= \int\limits_{\Omega} \varphi(y,\,|f(y)|)\,dy\quad \mbox{for all} \quad f\in L_0(\Omega).
$$

We put $L_{\varphi}= \left\{f\in L_0(\Omega):\; \rho_{\varphi}(\lambda_0 f)< \infty\quad \mbox{for some} \quad
\lambda_0> 0\right\}$ and
$$
\|f\|_{L_{\varphi}}= \inf\left\{\lambda> 0:\; \rho_{\varphi}\left(\frac f{\lambda}\right)\le 1\right\}.
$$
The space $L_{\varphi}$ is called Musielak-Orlicz space.
\end{definition}

Let $\omega$ be a weight function on $\Omega,$ i.e., $\omega$ is a non-negative, almost everywhere positive function on $\Omega.$
In this work we considered the weighted Musielak-Orlicz spaces.  We denote
\begin{equation*}
L_{\varphi,\,\omega}= \left\{f\in L_0(\Omega) :\,\;f\omega\in L_{\varphi}\right\}.
\end{equation*}
It is obvious that the norm in this spaces is given by
\begin{equation*}
\|f\|_{L_{\varphi,\, \omega}}= \|f\omega\|_{L_{\varphi}}.
\end{equation*}

\begin{remark}
Let $\varphi(x, t)= t^{q(x)}$ in the Definition 4, where $1\le q(x)< \infty$ and $x\in \Omega.$ Then we have the definition of variable exponent weighted Lebesgue spaces $L_{q(x)}\left(\Omega\right)$ (see [14]).
\end{remark}
{\bf Example 1.3.} The following Lemma shows that the Musielak-Orlicz spaces $L_{\varphi}$ is $p$-convex BFS.
\begin{lemma}{[4]} Let $\Omega_1\subset R^n$ and $\Omega_2\subset R^m$. Let $(x, t)\in \Omega_1\times [0, \infty),$ and $\varphi\left(x, t^{1/p}\right)\in \Phi$ for some $1\le p< \infty.$ Suppose $f: \Omega_1\times \Omega_2\mapsto R.$ Then the inequality
$$
\left\|\|f(x, \cdot)\|_{L_p\left(\Omega_2\right)}\right\|_{L_{\varphi}}\le
2^{1/p}\,\left\|\|f(\cdot, y)\|_{L_{\varphi}}\right\|_{L_p\left(\Omega_2\right)}
$$
is valid.
\end{lemma}

\begin{definition}{[38]} We say that $\varphi\in \Phi$ satisfies the $\Delta_2$-condition if there
exists $K \ge 2$ such that
$$
\varphi(y, 2t) \le K\,\varphi(y, t) \eqno (1.1)
$$
for all $y\in \Omega$ and all $t > 0.$ The smallest such $K$ is called the $\Delta_2$-constant
of $\varphi.$
\end{definition}

\begin{lemma} Let $\varphi \in\Phi$ and $1< s\le q(y)\le \overline q< \infty.$
Suppose for all $C> 0$ the condition
$$
\varphi(y, Ct) \le C^{q(y)}\,\varphi(y, t) \eqno (1.2)
$$
holds, where $y\in \Omega$ and $t > 0.$

Then a function $\varphi$ satisfies the $\Delta_2$-condition, with constant $K= 2^{\overline q}.$
\end{lemma}

{\it Proof.} Assume that (1.2) holds. Taking $C= 2$ in (1.2), we have
$$
\varphi(y, 2t) \le 2^{q(y)}\,\varphi(y, t)\le 2^{\overline q}\,\varphi(y, t).
$$
Thus the inequality (1.1) holds with constant $K= 2^{\overline q}.$

Lemma 3 is proved.

It is clear that if $\varphi(x, t)= t^{q(x)},$ then condition (1.2) satisfies automatically.

The following Lemma characterize bounded, sublinear operators from one Musielak-Orlicz spaces to another.
\begin{definition}
A mapping $\,S\,$ from one Musielak-Orlicz space $\,L_{\varphi}\,$ to another Musielak-Orlicz space $\,L_{\psi}\,$ is said to be sublinear if for all $f, g \in L_{\varphi}$  and $\lambda> 0$, we have

(1) $\;S (\lambda f) = \lambda\, S(f);$

(2) $\;S (f + g) \le S(f) + S(g).$
\end{definition}

\begin{lemma}{[15]} Let $\varphi, \psi\in \Phi$ and $\varphi$ and $\psi$ satisfy the $\Delta_2$ condition. Suppose  $S: L_{\varphi}\mapsto L_{\psi}$ be sublinear. Then the following conditions are equivalent:

(a) $\,$ $S$ is bounded, i.e. there exists $C> 0$ such that $\|Sf\|_{L_{\psi}}\le C\,\|f\|_{L_{\varphi}};$

(b) there exists $M_1, M_2> 0$ such that $\|f\|_{L_{\varphi}}\le C_1$ $\Longrightarrow$ $\rho_{L_{\psi}}\left(SF\right)\le M_2.$
\end{lemma}

We note that the Lebesgue spaces with mixed norm, weighted Lorentz spaces and e.t.c. is $p$-convex ($p$-concave) BFS.
Now we reduce more general result connected with Minkowski's integral inequality.

Let $X$ and $Y$ be BFSs on $\left(\Omega_1, \mu\right)$ and $\left(\Omega_2, \nu\right)$ respectively. By $X[Y]$ and $Y[X]$ we denote the spaces with mixed norm and consisting of all functions $g\in L_0\left(\Omega_1\times \Omega_2,\,\mu\times\nu\right)$ such that $\|g(x,\,\cdot)\|_Y\in X$ and $\|g(\cdot,\,y)\|_X\in Y.$ The norms in this spaces is defined as
$$
\|g\|_{X[Y]}= \left\|\|g(x,\,\cdot)\|_Y\right\|_X,\quad \|g\|_{Y[X]}= \left\|\|g(\cdot,\, y)\|_X\right\|_Y.
$$

\begin{theorem}{[43]} Let $X$ and $Y$ be BFSs with the Fatou property. Then the generalized Minkowski integral inequality
$$
\|f\|_{X[Y]}\le M\,\|f\|_{Y[X]}
$$
holds for all measurable functions $f(x,y)$ if and only if there exists $1\le p\le \infty$ such that $X$ is $p$-convex and $Y$ is $p$-concave.
\end{theorem}

It is known that $X[Y]$ and $Y[X]$ are BFSs on $\Omega_1\times \Omega_2$ (see [29].)

\begin{center}
3. {\bf Main results.}
\end{center}

We consider the multidimensional Hardy type operator and its dual operator
\begin{equation*}
Hf(x)= \int\limits_{|y|< |x|} f(y)\,dy \quad\mbox{and}\quad H^*f(x)=
\int\limits_{|y|> |x|} f(y)\,dy,
\end{equation*}
where $f\ge 0$ and $x\in R^n.$

Now we prove a two-weight criterion for multidimensional Hardy type
operator acting from the $p$-concave weighted BFS to weighted Lebesgue spaces.

\begin{theorem} Let $v(x)$ and $w(x)$ are weights on $R^n.$
Suppose that $X_w$ be a $p$-convex weighted BFSs for $1\le p< \infty$ on $R^n.$
Then the inequality
$$
\left\|Hf\right\|_{X_w}\le C
\,\left\|f\right\|_{L_{p,\, v}} \eqno(3.1)
$$
holds for every $f\ge 0$ and for all $\alpha\in (0, 1)$ if and only if
$$
A(\alpha)= \sup\limits_{t> 0} \left(\int\limits_{|y|< t}
[v(y)]^{- p'}\;dy\right)^{\frac{\alpha}{p'}}
\left\|\chi_{\{|z|> t\}}(\cdot)\left(\int\limits_{|y|< |\cdot|} [v(y)]^{-
p'}\;dy\right)^{\frac{1- \alpha}{p'}}\right\|_{X_w}< \infty. \eqno(3.2)
$$
Moreover, if $C> 0$ is the best possible constant in (3.1), then
\begin{equation*}
\sup\limits_{0< \alpha< 1} \frac{p'\,A(\alpha)}{(1- \alpha)\,\left[\left(\frac{
p'}{1-\alpha}\right)^{ p}+ \frac 1{\alpha\left(p- 1\right)}\right]^{1/p}} \le
C\le M\,\inf\limits_{0< \alpha< 1} \frac{A(\alpha)}{(1-
\alpha)^{1/p'}}.
\end{equation*}
\end{theorem}

{\it Proof of Theorem 2.} Sufficiency. Passing to the polar coordinates, we have
\begin{equation*}
h(y)= \left(\int\limits_{|z|<|y|} [v(z)]^{-
p'}\,dz\right)^{\frac \alpha{p'}}=\left(
\int\limits_0^{|y|} s^{n-1} \left(\int\limits_{|\xi|=
1}[v(s\xi)]^{- p' }\,d\xi\right) ds \right)^{\frac
\alpha{p'}},
\end{equation*}
where $d\xi$ is the surface element on the unit sphere. Obviously,
$h(y)= h(|y|),$ i.e., $h(y)$ is a radial function.

Applying H\"{o}lder's inequality for $L_{p}(R^n)$
spaces and after some standard transformations, we
have
\begin{equation*}
\begin{split}
\left\|Hf\right\|_{X_w}=&
\left\|w(\cdot)\,\int\limits_{|y|< |\cdot|} f(y)\, dy\right\|_{X}= \left\|w(\cdot)\,\int\limits_{|y|< |\cdot|} [f(y)h(y)v(y)] \left[h(y) v(y)\right]^{-1}dy \right\|_{X}\\
&\le \left\|w(\cdot)\,\left\|f\, h\,v\right\|_{L_{p}(|y|< |\cdot|)}
\left\|\left[h\,v\right]^{-1}\right\|_{L_{p'}( |y|<
|\cdot|)}\right\|_{X}\\
&= \left\|\left\|w(\cdot)\,f\, h\,v\,\chi_{\{|y|< |\cdot|\}}(y)
\left\|\left[h\,v \right]^{-1}\right\|_{L_{p'}(|y|<
|\cdot|)}\right\|_{L_{p}}\right\|_{X}\\
&= \left\|w\,f\, h\,v\,\chi_{\{|\cdot|< |x|\}}(\cdot)
\left\|\left[h\,v \right]^{-1}\right\|_{L_{p'}(|\cdot|<
|x|)}\right\|_{X\left[L_p\right]}.
\end{split}
\end{equation*}
Applying Theorem 1, we have
\begin{equation*}
\begin{split}
&\left\|w\,f\, h\,v\,\chi_{\{|\cdot|< |x|\}}(\cdot)
\left\|\left[h\,v \right]^{-1}\right\|_{L_{p'}(|\cdot|<
|x|)}\right\|_{X\left[L_p\right]}
\\
&\le M\,\left\|w\,f\, h\,v\,\chi_{\{|\cdot|< |x|\}}(\cdot)
\left\|\left[h\,v \right]^{-1}\right\|_{L_{p'}(|\cdot|<
|x|)}\right\|_{L_p\left[X\right]}
\\
&= M\,\left\|\left\|w(\cdot)\,f h\,v\, \chi_{\{|y|< |\cdot|\}}(y)
\left\|\left[h\,v\right]^{-1}\right\|_ {L_{p'}(|y|<
|\cdot|)}\right\|_{X}\right\|_{L_{p}}
\\
&= M\,\left\|f h\,v\,\left\|w(\cdot)\, \chi_{\{|y|< |\cdot|\}}(y)
\left\|\left[h\,v\right]^{-1}\right\|_ {L_{p'}(|y|<
|\cdot|)}\right\|_{X}\right\|_{L_{p}}.
\end{split}
\end{equation*}
By switching to polar coordinates and after some calculations, we
get
$$
\left\|\left[h\,v\right]^{-1}\right\|_{L_{p'}\left(|y|< |x|\right)} =
\left(\int\limits_{|y|< |x|} [h(|y|)\,v(y)]^{-p'}\, dy\right)^{1/p'}
$$
$$
= \left(\int\limits_0^{|x|} r^{n-1}\, [h(r)]^{-p'}\,\left[\int\limits_{|\xi|= 1}
\left[v(r\xi)\right]^{- p'}\,d\xi\right]\,dr\right)^{1/p'}
$$
$$
= \left(\int\limits_0^{|x|} \left[\int\limits_0^r s^{n- 1}
\left(\int\limits_{|\xi|= 1}[v(s\xi)]^{-p'}
d\xi\right) ds\right]^{-\alpha} \left(\int\limits_{|\xi|= 1}
\left[v(r\xi)\right]^{-p'}\, d\xi\right) r^{n-1}
dr\right)^{1/p'}
$$
$$
= \frac 1{\left(1- \alpha\right)^{1/p'}}\,
\left(\int\limits_0^{|x|} \frac{d}{dr}\left\{\left(\int\limits_0^r
s^{n- 1} \left(\int\limits_{|\xi|= 1} [v(s\xi)]^{-p'}
d\xi\right) ds\right)^{1- \alpha}\right\} dr\right)^{1/
p'}
$$
$$
= \frac 1{\left(1- \alpha\right)^{1/p'}}\,
\left(\int\limits_0^{|x|} s^{n- 1}\,\left(\int\limits_{|\xi|= 1}
[v(s\xi)]^{-p'} d\xi\right) ds\right)^{\frac {1-
\alpha}{p'}}= \frac 1{\left(1- \alpha\right)^{1/p'}}\,
\left(\int\limits_{|z|< |x|} [v(z)]^{-
p'}\,dz\right)^{\frac {1- \alpha}{p'}}.
$$
Therefore from the condition (3.2), we obtain
$$
\left\|f h\,v\,\left\|w(\cdot)\, \chi_{\{|y|< |\cdot|\}}(y)
\left\|\left[h\,v\right]^{-1}\right\|_ {L_{p'}(|y|<
|\cdot|)}\right\|_{X}\right\|_{L_{p}}
$$
$$
= \frac 1{\left(1- \alpha\right)^{1/p'}}\,\left\|f\,v\,\left[h\left\|\chi_{\{|\cdot|> |y|\}} \,
\left(\int\limits_{|z|< |\cdot|} [v(z)]^{-p'}\,dz\right)^{\frac {1- \alpha}{p'}}\right\|_
{X_w}\right]\right\|_{L_{p}}\le \frac {A(\alpha)}{\left(1- \alpha\right)^{1/
p'}}\, \left\|f\,v\right\|_{L_{p}}.
$$
Thus
\begin{equation*}
\left\|Hf\right\|_{X_w}\le
M\,\frac {A(\alpha)}{\left(1- \alpha\right)^{1/
p'}}\, \left\|f\right\|_{L_{p,\,v}}\; \mbox{for all}\; \alpha\in (0, 1).
\end{equation*}

{\bf Necessity.} Let $f\in L_{p, v}\left(R^n\right),$ $f\ge 0$
and the inequality (3.1) is valid. We choose the test function as
\begin{equation*}
f(x)= \frac{p'}{1- \alpha}\,[g(t)]^{-\frac
{\alpha}{p'}- \frac 1{p}}\,v^{-
p'}(x)\,\chi_{\{|x|< t\}}(x)+ [g(|x|)]^{-\frac {\alpha}{
p'}- \frac 1{p}} v^{-p'}(x)\,\chi_{\{|x|>
t\}}(x),
\end{equation*}
where $t> 0$ is a fixed number and
\begin{equation*}
g(t)= \int\limits_{|y|< t} v^{-p'}(y)\, dy=
\int\limits_0^t s^{n- 1}\left(\int\limits_{|\eta|= 1} v^{-
p'}(s \eta)\, d\eta\right)\,ds.
\end{equation*}
It is obvious that $\displaystyle{\frac{dg}{dt}= t^{n-
1}\,\int\limits_{|\eta|= 1} v^{-p'}(t \eta)\, d\eta}.$
Again by switching to polar coordinates, from the right hand side of
inequality (3.1) we get that
$$
\left\|f\right\|_{L_{p, v}}=\left[ \int\limits_{|x|<
t}\left(\frac{p'}{1- \alpha}\right)^{p}
[g(t)]^{-\alpha(p- 1) - 1}\,v^{-p'}(x)\, dx +
\int\limits_{|x|> t}[g(|x|)]^{-\alpha(p- 1) -
1}\,v^{-p'}(x)\, dx\right]^ {1/p}
$$
$$
= \left[\left(\frac{p'}{1- \alpha}\right)^{
p}\,[g(t)]^{\alpha(1- p)}\, + \int\limits_{t}^{\infty}
r^{n- 1}\,[g(r)]^{- \alpha(p- 1) -
1}\,\left(\int\limits_{|\xi|= 1} v^{-
p'}(r\xi)\,d\xi\right)\,dr\right]^{1/ p}
$$
$$
= \left[\left(\frac{ p'}{1- \alpha}\right)^{
p}\,[g(t)]^{\alpha(1-  p)}- \frac{1 }{\alpha( p-
1)} \int\limits_{t}^{\infty} \frac d{dr}[g(r)]^{-\alpha(
p- 1)}\,dr \right]^{1/ p}
$$
$$
=\left[\left(\frac{ p'}{1- \alpha}\right)^{
p}\,[g(t)]^{\alpha(1-  p)}+ \frac{1 }{\alpha( p-
1)}\left\{ [g(t)]^{-\alpha( p- 1)}-
\left[\int\limits_{R^n} v^{- p'}(y)\,
dy\right]^{-\alpha( p- 1)}\, \right\}\right]^{1/ p}
$$
$$
\le \left[\left(\frac{ p'}{1- \alpha}\right)^{
p}+ \frac 1{\alpha( p- 1)}\right]^{1/ p}
[g(t)]^{-\frac{\alpha}{ p'}}= \left[\left(\frac{
p'}{1- \alpha}\right)^{ p}+ \frac 1{\alpha( p-
1)}\right]^{1/ p}\, [h(t)]^{-1}.
$$
After some calculations, from the left hand side of inequality (3.1), we have
\begin{equation*}
\left\|Hf \right\|_{X_w}= \left\|\int\limits_{|y|< |\cdot|} f(y)\, dy\right\|_{X_w}
 \ge \left\|\chi_{\{|\cdot|> t\}}\,\int\limits_{|y|< |\cdot|} f(y)\,dy\right\|_{X_w}=
\end{equation*}
\begin{equation*}
\begin{split}
&= \left\| \chi_{\{|\cdot|> t\}}\left(\frac{p'}{1- \alpha} \,\int\limits_{|y|< t}
[g(t)]^{-\frac {\alpha}{p'}- \frac 1{p}}
v^{-p'}(y)\, dy + \int\limits_{t< |y|< |\cdot|}
[g(|y|)]^{-\frac {\alpha}{p'}- \frac 1{p}}
v^{-p'}(y)\, dy\right) \right\|_{X_w}\\
&= \left\| \chi_{\{|\cdot|> t\}}\left(\frac{p'}{1- \alpha}[g(t)]^{\frac {1-
\alpha}{p'}}+ \int\limits_{t}^{|\cdot|} r^{n-
1}\,[g(r)]^{-\frac{\alpha}{p'}- \frac 1{
p}} \left(\int\limits_{|\eta|= 1} v^{-p'}(r\eta)\,d\eta
\right)\,dr\right)\right\|_{X_w}\\
&= \left\|\chi_{\{|\cdot|> t\}}\left(\frac{p'}{1- \alpha}[g(t)]^{\frac {1-
\alpha}{p' }}+ \frac{p'}{1- \alpha}\,
\int\limits_{t}^{|\cdot|} \frac d{dr} [g(r)]^{\frac {1-
\alpha}{p'}}\,dr\right)\right\|_{X_w}\\
&\left\| \chi_{\{|\cdot|> t\}}\,\left[\frac{p'}{1- \alpha}[g(t)]^{\frac {1-
\alpha}{p'}}+ \frac{p'}{1- \alpha}\,
\left([g(|\cdot|)]^{\frac {1- \alpha}{p'}}- [g(t)]^{\frac
{1- \alpha}{p'}}\right)\right]\right\|_{X_w}\\
&=\frac{p'}{1- \alpha}\,\left\|\chi_{\{|\cdot|> t\}}\, [g(\cdot)]^{\frac
{1- \alpha}{p'}} \right\|_{X_w}.
\end{split}
\end{equation*}
Hence, this implies that
\begin{equation*}
\frac{p'}{1- \alpha}\,\left[\left(\frac{p'}{1-
\alpha}\right)^{p}+ \frac 1{\alpha(p-
1)}\right]^{-1/ p}\, [g(t)]^{\frac{\alpha}{
p'}}\,\left\|\chi_{\{|\cdot|> t\}}\,[g(\cdot)]^{\frac {1- \alpha}{
p'}}\right\|_{X_w} \le C,
\end{equation*}
i.e., $\displaystyle{\frac{p'\,A(\alpha)}{(1-
\alpha)\,\left[\left(\frac{p'}{1-\alpha}\right)^{p}+ \frac
1{\alpha\left(p- 1\right)}\right]^{1/p}}\le
C}$ for all $\alpha\in (0, 1).$

This completes the proof of Theorem 2.

For the dual operator, the below stated  theorem is proved analogously.

\begin{theorem} Let $v(x)$ and $w(x)$ are weights on $R^n.$
Suppose that $X_w$ be a $p$-convex weighted BFSs for $1\le p< \infty$ on $R^n.$
 Then the inequality
$$
\left\|H^*f\right\|_{X_w}\le C
\,\left\|f\right\|_{L_{p,\, v}} \eqno(3.3)
$$
holds for every $f\ge 0$ and for all $\gamma\in (0,
1)$ if and only if
\begin{equation*}
B(\gamma)= \sup\limits_{t> 0} \left(\int\limits_{|y|> t}
[v(y)]^{- p'}\;dy\right)^{\frac{\gamma}{p'}}
\left\|\chi_{\{|z|< t\}}(\cdot) \left(\int\limits_{|y|> |\cdot|} [v(y)]^{-
p'}\;dy\right)^{\frac{1- \gamma}{p'}}\right\|_{X_w}< \infty.
\end{equation*}
Moreover, if $C> 0$ is the best possible constant in (3.3) then
\begin{equation*}
\sup\limits_{0< \gamma< 1} \frac{p'\,B(\gamma)}{(1- \gamma)\,\left[\left(\frac{
p'}{1- \gamma}\right)^{ p}+ \frac 1{\gamma\left(p- 1\right)}\right]^{1/p}} \le
C\le M\,\inf\limits_{0< \gamma< 1} \frac{B(\gamma)}{(1- \gamma)^{1/p'}}.
\end{equation*}
\end{theorem}

\begin{corollary}
Note that Theorem 2 and Theorem 3 in the case $X_w= L_{\varphi,\,w},$ $\varphi\left(x, t^{1/p}\right)\in \Phi$ for some $1\le p< \infty,$ $x\in R^n$ was proved in [4]. In the case $X_w= L_{q,\,w},$ $1< p\le q< \infty,$ for $x\in (0, \infty),$ $\displaystyle{\alpha= \frac{s- 1}{p- 1}}$ and $s\in (1,\,p)$ Theorem 2 and Theorem 3 was proved in [46]. For $x\in R^n$ in the case $X_w= L_{q(x),\,w}$ and $1< p\le q(x)\le ess\sup\limits_{x\in R^n} q(x)< \infty$ Theorem 2 and Theorem 3 was proved in [3] (see also [2]).
\end{corollary}

\begin{remark}  In the case $n=1,$ $X_w= L_{q,\,w},$ $1< p\le q\le \infty,$ at
$x\in (0, \infty),$ for classical Lebesgue spaces the various variants of Theorem 2 and
Theorem 3 were proved in [20], [10], [27], [28], [36], [37], [45] and etc. In particular,
in the Lebesgue spaces with variable exponent the boundedness of Hardy type operator was proved in
[13], [14], [17], [21], [26], [34], [35] and etc. For $X_w= L_{q(x),\,w},$ $1< p\le q(x)\le ess\sup\limits_{x\in [0, 1]} q(x)< \infty$ and $x\in [0, 1]$ the two-weighted criterion for one-dimensional Hardy operator was proved in [26]. Also, other type two-weighted criterion for multidimensional Hardy type operator in the case $X_w= L_{q(x),\,w},$ $1< p\le q(x)\le ess\sup\limits_{x\in R^n} q(x)< \infty$ and $x\in R^n$ was proved in [34] (see also [35]).  In the papers [9] and [42] the inequalities of modular type for more general operators was proved. Also, in [11] the Hardy type inequalities with special power-type weights in Orlicz spaces was proved.
\end{remark}

\vspace{5mm}
\begin{center}
4. {\bf Applications.}
\end{center}

Now we consider the multidimensional geometric mean operator defined as
\begin{equation*}
Gf(x)= \exp\left(\frac{1}{|B(0,\,|x|)|}\int\limits_{B(0,\,|x|)} \ln\,f(y)\,dy\right),
\end{equation*}
where $f> 0$ and $\displaystyle{|B(0,\,|x|)|= |B(0, 1)|\,|x|^n}.$
It is obvious that $G\left(f_1\cdot f_2\right)(x)= G f_1(x)\cdot Gf_2(x).$

We formulate a two-weighted criterion on boundedness of multidimensional geometric mean operator in weighted Musielak-Orlicz spaces.

\begin{theorem} Let $\varphi\left(x, t^{1/p}\right)\in \Phi$ for some $0< p< \infty$ and $x\in R^n.$ Suppose that
$v(x)$ and $w(x)$ are weight functions on $R^n.$ Then the inequality
$$
\left\|Gf\right\|_{L_{\varphi,\, w}}\le C
\,\left\|f\right\|_{L_{p,\,v}} \eqno(4.1)
$$
holds for every $f>0$ and for all $s\in (1, p)$ if and only if
$$
D(s)= \!\sup\limits_{t> 0} |B(0, t)|^{\frac{s- 1} p}\left\|\frac{\chi_{\{|z|> t\}}(\cdot)}{|B(0, |\cdot|)|^{\frac {s}p}}\,\exp\left(\frac 1{|B(0, |\cdot|)|}\int\limits_{B(0,|\cdot|)} \ln\frac 1{v(y)}\,dy\right)\right\|_
{L_{\varphi,\,w}}< \infty. \eqno(4.2)
$$
Moreover, if $C> 0$ is the best possible constant in (4.1) then
\begin{equation*}
\sup\limits_{s> 1} \frac{e^{\frac sp}}{\left(e^s+ \frac 1{s- 1}\right)^{1/p}}\,
D(s) \le C \le 2^{1/p}\,\inf\limits_{s> 1} e^{\frac{s- 1}{p}}\, D(s).
\end{equation*}
\end{theorem}

{\it Proof of Theorem 4.}
Let $\displaystyle{\alpha= \frac{s- 1}{p- 1}},$ where $1< s< p.$  We replace  $f$ with $f^{\beta},$ $v$ with $v^{\beta},$ $w$ with $\displaystyle{\frac {w^{\beta}(x)}{|B(0, |x|)|}},$ $0< \beta< p,$ and  $\,p\,$ with $\displaystyle{\frac p{\beta}}$ and $\varphi(x,\,t)\,$ with $\varphi\left(x,\,t^{1/\beta}\right)$ in (3.1), (3.2), we find that for $\displaystyle{1< s< \frac p{\beta}}$
\begin{equation*}
\begin{split}
\left\|\frac{w^{\beta}}{|B(0, |\cdot|)|}\,H(f^{\beta})\right\|_{L_{\varphi\left(\cdot, t^{1/\beta}\right)}}&= \left\|\left(\frac 1{|B(0, |\cdot|)|}
\int\limits_{B(0,\,|\cdot|)} f^{\beta}(y)\,dy\right)^{1/\beta}\right\|^{\beta}_{L_{\varphi,\, w}\left(R^n\right)}\\
&\le C_{\beta}\left(\int\limits_{R^n} [f(y) v(y)]^{p}\,dy\right)^{\beta/p}.
\end{split}
\end{equation*}
Then the inequality
$$
\left\|\left(\frac 1{|B(0, |\cdot|)|} \int\limits_{B(0,\,|\cdot|)} f^{\beta}(y)\,dy\right)^{1/\beta}
\right\|_{L_{\varphi,\, w}\left(R^n\right)}\le C_{\beta}^{1/\beta}\left(\int\limits_{R^n} [f(y) v(y)]^{p}\,dy\right)^{1/p} \eqno(4.3)
$$
holds if and only if
\begin{equation*}
\begin{split}
&A\left(\frac{s- 1}{p- 1}\right)\\
&
=\left[\sup\limits_{t> 0} \left(\int\limits_{|y|< t} [v(y)]^{-\frac{\beta\,p}{p- \beta}}\,dy\right)^
{\frac{s- 1}p}\left\|\left(\frac {\chi_{\{|z|> t\}}(\cdot)}{|B(0, |\cdot|)|^{\frac{p}{p- \beta s}}}
\int\limits_{|y|< |\cdot|} [v(y)]^{-\frac{\beta p}{p- \beta}}\,dy\right)^
{\frac{p- \beta s}{\beta p}} \right\|_{L_{\varphi,\, w}}\right]^{\beta}\\
&= B^{\beta}(s,\,\beta)< \infty
\end{split}
\end{equation*}
and
$$
\sup\limits_{1< s< \frac p\beta} \left[\frac{\left(\frac{p}{p- s\beta}\right)^{\frac p\beta}}
{\left(\frac{p}{p- s\beta}\right)^{\frac p\beta}+ \frac 1{s- 1}}\right]^{\beta/p} B^{\beta}(s, \beta)\le C_{\beta}
\le 2^{\frac {\beta}{p}}\,\inf\limits_{1< s< \frac p\beta} \left(\frac{p- \beta}{p- s\beta}\right)^{\frac {p- \beta}{p}}
\,B^{\beta}\left(s, \beta\right). \eqno(4.4)
$$
By the L'Hospital rule, we get
\begin{equation*}
\begin{split}
&\lim\limits_{\beta\to +0} \left(\frac 1{|B(0, |x|)|^{\frac{p}{p- \beta s}}}\int\limits_{|y|< |x|} [v(y)]^{-\frac{\beta p}
{p- \beta}}\,dy\right)^{\frac{p- \beta s}{\beta p}}\\
&= \lim\limits_{\beta\to +0} \exp\left[\frac{p\,\ln\frac 1{|B(0, |x|)|}
+ (p- \beta s)\,\ln\left(\int\limits_{|y|< |x|}[v(y)]^{-\frac{\beta p}{p- \beta}}\,dy\right)}{p\,\beta}\right]\\
&= \lim\limits_{\beta\to +0} \exp\left[-\frac sp\,\ln\left(\int\limits_{|y|< |x|}[v(y)]^{-\frac{\beta p}{p- \beta}}\,
dy\right)+ \frac{(p- \beta s)\left(\frac{p}{p- \beta}\right)^2\int\limits_{|y|< |x|} [v(y)]^{-\frac{\beta p}{p- \beta}}
\ln \frac 1{v(y)}\,dy}{p\,\int\limits_{|y|< |x|}[v(y)]^{-\frac{\beta p}{p- \beta}}\,dy}\right]\\
&= \exp\left[\frac sp\,\ln\frac 1{|B(0, |x|)|}+ \frac{\int\limits_{|y|< |x|} \ln \frac 1{v(y)}\,dy}{|B(0, |x|)|}\right]=
\frac 1{|B(0, |x|)|^{\frac sp}}\,\exp\left(\frac 1{|B(0, |x|)|}\int\limits_{B(0,|x|)} \ln\frac 1{v(y)}\,dy\right).
\end{split}
\end{equation*}
Therefore
\begin{equation*}
\begin{split}
\lim\limits_{\beta\to +0} B\left(s, \beta\right)&= \sup\limits_{t> 0} |B(0, t)|^{\frac{s- 1} p} \left\|\frac{\chi_{\{|z|> t\}}(\cdot)}{|B(0, |\cdot|)|^{\frac {s}p}}
\,\exp\left(\frac 1{|B(0, |\cdot|)|}\int\limits_{B(0,|\cdot|)} \ln\frac 1{v(y)}\,dy\right)\right\|_
{L_{\varphi,\,w}}\\
&= D(s)< \infty
\end{split}
\end{equation*}
and
$$
\sup\limits_{s> 1} \frac{e^{\frac sp}}{\left(e^s+ \frac 1{s- 1}\right)^{1/p}}\,D(s)\le
\lim\limits_{\beta\to +0} C_{\beta}^{1/\beta}\le 2^{1/p}\,\inf\limits_{s> 1} e^{\frac {s- 1}{p}}\,D(s). \eqno(4.5)
$$
Further, we have
\begin{equation*}
\lim\limits_{\beta\to +0} \left(\frac 1{|B(0, |x|)|}\int\limits_{B(0,\,|x|)} f^{\beta}(y)\,dy\right)^{1/\beta}=
\exp\left(\frac{1}{|B(0,\,|x|)|}\int\limits_{B(0,\,|x|)} \ln\,f(y)\,dy\right)= Gf(x).
\end{equation*}
From (4.4) it follows that $\displaystyle{\lim\limits_{\beta\to +0} C_{\beta}= 1},$ and according to (4.2) and (4.5)
$\displaystyle{\lim\limits_{\beta\to +0} C_{\beta}^{1/\beta}= C< \infty}.$ Therefore the inequality (4.1) is valid.
Moreover, from (4.3) for $\beta\to +0$ we obtain that
\begin{equation*}
\left\|Gf\right\|_{L_{q(\cdot),\, w}\left(R^n\right)}\le C
\,\left\|f\right\|_{L_{p,\,v}(R^n)}
\end{equation*}
and by (4.5)
\begin{equation*}
\sup\limits_{s> 1} \frac{e^{\frac sp}}{\left(e^s+ \frac 1{s- 1}\right)^{1/p}}\,
D(s) \le C \le 2^{1/p}\,\inf\limits_{s> 1} e^{\frac{s- 1}{p}}\, D(s).
\end{equation*}

This completes the proof of Theorem 4.

\begin{remark}
Let $\varphi(x, t)= t^q$ and $n= 1.$ Note that the simplest case of (2.4) with $v= w= 1$ and $p= q= 1$
was considered in [20] and in [25]. Later this inequality was generalized in various ways by many authors in
[12, 22, 23, 24, 31, 39, 40, 41, 46] and etc.
\end{remark}

\begin{corollary}
Let $\varphi(x, t)= t^q,$ $0< p\le q< \infty$ and let $f$ be a positive function on $R^n.$ Then
$$
\left(\int\limits_{R^n} [Gf(x)]^q\,|x|^{\delta\,q}\;dx\right)^{1/q}\le
C\;\left(\int\limits_{R^n} f^p (x)\,|x|^{\mu\,p}\;dx\right)^{1/p} \eqno(4.6)
$$
holds with a finite constant $C$ if and only if
\begin{equation*}
\delta+ \frac nq= \frac{\mu}n+ \frac np
\end{equation*}
and the best constant $C$ has the following condition:
\begin{equation*}
\sqrt[q]{\frac p{nq}}\;e^{\frac{\mu}{n^2}}\;|B(0,\,1)|^{\frac 1q- \frac 1p}\,\sup\limits_{s> 1}
\frac{e^{\frac sp}\,(s- 1)^{\frac 1p- \frac 1q}}{\left[(s- 1)e^s+ 1\right]^{1/p}}\le C\le
\frac{|B(0,\,1)|^{\frac 1q- \frac 1p}\;e^{\frac{\mu}{n^2}+ \frac 1q}}{\sqrt [q]{n}}.
\end{equation*}
\end{corollary}

\begin{remark}
Let $\varphi(x, t)= t^q$ and $q= p.$ Then the inequality (4.6) is sharp with the constant
$\displaystyle{C= \frac {e^{\frac{\mu}{n^2}+ \frac 1p}} {\sqrt [p]{n}}}.$
\end{remark}

The sufficient conditions for general weights ensuring the validity
of the two-weight strong type inequalities for some sublinear
operator are given in the following theorem.

\begin{theorem} Let $\varphi\left(x, t^{1/p}\right)\in \Phi$ for some $1< p< \infty$ and $x\in R^n.$ Suppose that
$v(x)$ and $w(x)$ are weight functions on $R^n.$ Let $T$ be a sublinear operator acting boundedly from $L_{p}\left(R^n\right)$ to $L_{\varphi}\left(R^n\right)$ such that, for any $f\in L_1(R^n)$ with compact support and $x\notin supp\; f$
$$
|Tf(x)|\le C\,\int\limits_{R^n} \frac{|f(y)|}{|x- y|^n}\,dy,
\eqno(4.7)
$$
where $C> 0$ is independent of $f$ and $x.$ Let there exists $1< p\le r(x)\le ess\,\sup\limits_{x\in R^n} r(x)<\infty$ such that, for all $C> 0$ $\varphi \left(y, Ct\right)\le C^{r(y)}\,\varphi (y, t).$

Moreover, let $v(x)$ and $w(x)$ are weight functions on $R^n$ and
satisfies the following conditions:
$$
A= \sup\limits_{t> 0}\left(\int\limits_{|y|< t} [v(y)]^{-
 p^{\,'}}\;dy\right)^{\frac{\alpha}{p{\,'}}}
\left\|\frac{\chi_{\{|x|> t\}}}{|x|^n}\left(\int\limits_{|y|< |x|}
[v(y)]^{- p^{\,'}}\;dy\right)^{\frac{1- \alpha}{p^{\,'}}}\right\|_{L_{\varphi,\,w}}< \infty, \eqno(4.8)
$$
$$
B= \sup\limits_{t> 0}\left(\int\limits_{|y|> t} [v(y)|y|^n]^{-
p^{\,'}}\;dy\right)^{\frac{\beta}{p{\,'}}}
\left\|\chi_{\{|x|< t\}}\left(\int\limits_{|y|> |x|} [v(y)|y|^n]^{-
p^{\,'}}\;dy\right)^{\frac{1- \beta}{p^{\,'}}}\right\|_{L_{\varphi,\,w}}< \infty. \eqno(4.9)
$$

There exists $M> 0$ such that
$$
\sup\limits_{|x|/2< |y|\le 4\,|x|} w(y)\le M\, \inf\limits_{|x|/2< |y|\le 4\,|x|}
v(x). \eqno(4.10)
$$
Then there exists a positive constant $C,$ independent of $f$, such
that for all $f\in L_{p, v}(R^n)$
$$
\|Tf\|_{L_{\varphi, w}}\le C \|f\|_{L_{p, v}(R^n)}.
$$
\end{theorem}

{\it Proof of Theorem 5.} Let $Z= \left\{0,\pm 1,\pm 2,\ldots,\right\}.$ For $k\in Z$ we define $E_k= \left\{x\in R^n:\right.$ $\left.2^k< |x|\le 2^{k+ 1}\right\},$ $E_{k,1}= \left\{x\in R^n:\;|x|\le 2^{k- 1}\right\},$ $E_{k,2}= \left\{x\in
R^n:\;2^{k- 1}< |x|\le 2^{k+ 2}\right\},$ $E_{k, 3}= \left\{x\in R^n:\; |x|> 2^{k- 1}\right\}.$ Then $E_{k, 2}= E_{k- 1}\cup E_k\cup E_{k+ 1}$ and the multiplicity of the covering $\left\{E_{k,
2}\right\}_{k\in Z}$ is equal to 3.

Given $f\in L_{p, v}(R^n),$ we write
$$
|Tf(x)|= \sum\limits_{k\in Z} |Tf(x)|\,\chi_{E_k}(x)\le
\sum\limits_{k\in Z} \left|Tf_{k,1}(x)\right|\,\chi_{E_k}(x)+
\sum\limits_{k\in Z} \left|Tf_{k,2}(x)\right|\,\chi_{E_k}(x)+
$$
$$
+ \sum\limits_{k\in Z} \left|Tf_{k,3}(x)\right|\,\chi_{E_k}(x) = T_1
f(x)+ T_2 f(x)+ T_3 f(x),
$$
where $\chi_{E_k}$ is the characteristic function of the set $E_k,$
$f_{k,i}= f\chi_{E_{k,i}},$ $i= 1, 2, 3.$

First we shall estimate $\left\|T_1 f\right\|_{L_{p, \omega_2}}.$
Note that for $x\in E_k,$ $y\in E_{k, 1}$ we have $|y|< 2^{k- 1}\le
|x|/2.$ Moreover, $E_k\cap supp\,f_{k,1}= \emptyset$ and $|x- y|\ge |x|- |y|\ge
|x|- |x|/2= |x|/2.$ Hence by (4.7)
$$
\left|T_1 f(x)\right|\le C\,\sum\limits_{k\in
Z}\left(\int\limits_{R^n} \frac{|f_{k, 1}(y)|}{|x-
y|^n}\,dy\right)\,\chi_{E_k}\le C\,\int\limits_{|y|< |x|/2}
\frac{|f(y)|}{|x- y|^n}\,dy\le
$$
$$
\le C\,\int\limits_{|y|< |x|} \frac{|f(y)|}{|x- y|^n}\,dy \le
2^n\,C\,|x|^{-n}\,\int\limits_{|y|< |x|} |f(y)|\,dy
$$
for any $x\in E_k.$ Hence we have
$$
\left\|T_1 f\right\|_{L_{\varphi,\,w}}\le 2^n\,C\,\left\||x|^{-n}\,\int\limits_{|y|< |x|} |f(y)|\,dy
\right\|_{L_{\varphi,\,w}}= \left\|\int\limits_{|y|< |x|} |f(y)|\,dy \right\|_{L_{\varphi,\,|x|^{-n}\,w}}.
$$
By the condition (4.8) and Theorem 2, we obtain
$$
\left\|T_1 f\right\|_{L_{\varphi,\,w}}\le C_1\,\|f\|_{L_{p,
v}(R^n)} \eqno(4.11)
$$
where $C_1> 0$ is independent of $f$ and $x\in R^n.$

Next we estimate $\left\|T_3 f\right\|_{L_{p, w}(R^n)}.$ It is
obviously that, for $x\in E_k,$ $y\in E_{k, 3}$ we have $|y|>
2\,|x|$ and $|x- y|\ge |y|- |x|\ge |y|- |y|/2= |y|/2.$ Since
$E_k\cap\,supp\,f_{k,3}= \emptyset,$ for $x\in E_k$ by (4.7), we have
$$
\left|T_3 f(x)\right|\le C\,\int\limits_{|y|> 2|x|}
\frac{|f(y)|}{|x- y|^n}\,dy\le 2^n\,C\,\int\limits_{|y|> 2|x|}
\frac{|f(y)|}{|y|^n}\,dy.
$$
Hence we obtain
$$
\left\|T_3 f\right\|_{L_{\varphi,\,w}}\le 2^n\,C\,
\left\|\int\limits_{|y|> 2|\cdot|} |f(y)|\,|y|^{-n}\,dy
\right\|_{L_{\varphi,\,w}}\le
$$
$$
\le 2^n\,C\, \left\|\int\limits_{|y|> |x|} |f(y)|\,|y|^{-n}\,dy
\right\|_{L_{\varphi,\,w}}.
$$
By the condition (4.9) and Theorem 3, we obtain
$$
\left\|T_3 f\right\|_{L_{\varphi,\,w}}\le C_2\,\|f\|_{L_{p,v}(R^n)},\eqno(4.12)
$$
where $C_2> 0$ is independent of $f$ and $x\in R^n.$

Let $T:L_{p}\left(R^n\right)\to L_{\varphi}\left(R^n\right).$  We have
$$
\left\|Tf_{k, 2}\right\|_{L_{\varphi,\,w}\left(R^n\right)}= \left\|\sum\limits_{k\in Z} \left|Tf_{k, 2}\right|\,\chi_{E_k} \right\|_{L_{\varphi,\,w}\left(R^n\right)}.
$$
By virtue of Lemma 3 it suffices to prove that from $\|f\|_{L_{p,\,v}\left(R^n\right)}\le 1$ implies \linebreak $\displaystyle{\int\limits_{R^n} \varphi \left(y, w\,\sum\limits_{k\in Z} \left|Tf_{k, 2}\right|\,\chi_{E_k}\right) \,dx\le C},$ where $C> 0$ is independent on $k\in Z.$

Finally, we estimate $\left\|T_2 f\right\|_{L_{\varphi,\,w}}.$ By
the $L_{p}(R^n)\mapsto L_{\varphi}\left(R^n\right)$ boundedness of $T$ and condition (4.10), we have
$$
\int\limits_{R^n} \varphi \left(y, w(y)\,\sum\limits_{k\in Z} \left|Tf_{k, 2}(y)\right|\,\chi_{E_k}(y)\right) \,dy=
\sum\limits_{m\in Z} \int\limits_{E_m} \varphi \left(y, w(y)\,\sum\limits_{k\in Z} \left|Tf_{k, 2}(y)\right|\,\chi_{E_k}(y)\right) \,dy=
$$
$$
\sum\limits_{k\in Z} \int\limits_{E_k} \varphi \left(y, w(y)\,\left|Tf_{k, 2}(y)\right|\right) \,dy=
\sum\limits_{k\in Z} \int\limits_{E_k} \varphi \left(y, C\,w(y) \left\|f_{k, 2}\right\|_{L_p\left(R^n\right)} \,\frac{\left|Tf_{k, 2}\right|}{C\left\|f_{k, 2}\right\|_{L_p\left(R^n\right)}}\right) \,dy\le
$$
$$
C_1\,\sum\limits_{k\in Z}\, \int\limits_{E_k} \left(C\,w(y) \left\|f_{k, 2}\right\|_{L_p\left(R^n\right)}\right)^{r(y)} \varphi \left(y, \frac{\left|Tf_{k, 2}\right|}{C\left\|f_{k, 2}\right\|_{L_p\left(R^n\right)}}\right) \,dy\le
$$
$$
C_2\,\sum\limits_{k\in Z}\,\sup\limits_{y\in E_k} \left(w(y) \left\|f\right\|_{L_p\left(E_{k,2}\right)} \right)^{r(y)}\,\int\limits_{R^n} \varphi \left(y, \frac{\left|Tf_{k, 2}\right|}{C\left\|f_{k,2}\right\|_ {L_p\left(R^n\right)}}\right) \,dy\le
$$
$$
C_2\,\sum\limits_{k\in Z}\,\sup\limits_{y\in E_k} \left(w(y) \left\|f\right\|_{L_p\left(E_{k,2}\right)} \right)^{r(y)}= C_2\,\sum\limits_{k\in Z}\,\sup\limits_{y\in E_k} \left(\|f\,w\|_{L_{p}\left(E_{k,2}\right)} \right)^{r(y)}\le
$$
$$
C_3 \sum\limits_{k\in Z}\,\sup\limits_{y\in E_k} \left(\|f\,\inf\limits_{y\in E_{k,2}}\,v(y)\|_{L_{p}\left(E_{k,2}\right)}\right)^{r(y)}
\le C_3 \sum\limits_{k\in Z}\,\sup\limits_{y\in E_k} \left(\|f\,v\|_{L_{p}\left(E_{k,2}\right)}\right)^{r(y)}=
$$
$$
C_3 \sum\limits_{k\in Z} \left(\|f\|_{L_{p,\,v}\left(E_{k,2}\right)}\right)^{\inf\limits_{y\in E_k} r(y)}\le
C_3 \sum\limits_{k\in Z} \left(\|f\|_{L_{p,\,v}\left(E_{k,2}\right)}\right)^{\underline r}=
$$
$$
= C_3 \sum\limits_{k\in Z} \left(\int\limits_{E_{k, 2}} \left[|f(y)|\,v(y)\right]^p\, dy\right)^{\underline r/p}\le
C_3 \left(\sum\limits_{k\in Z} \int\limits_{E_{k, 2}} \left[|f(y)|\,v(y)\right]^p\, dy\right)^{\underline r/p}=
$$
$$
C_3 \,\left(\sum\limits_{k\in Z} \left[\int\limits_{E_{k- 1}}+ \int\limits_{E_{k}}+ \int\limits_{E_{k+1}}\right] \left[|f(y)|\,v(y)\right]^p\,dy\right)^{\underline r/p}= C_3\,\left(3\,\sum\limits_{k\in Z} \int\limits_{E_k}
\left[|f(y)|\,v(y)\right]^p\,dy\right)^{\underline r/p}=
$$
$$
3^{{\underline r}/p}\, \left(\|f\|_{L_{p,\,v}\left(R^n\right)}^p\right)^{\underline r/p}=
C_4\,\|f\|_{L_{p,\,v}\left(R^n\right)}^{\underline r}\le C_4.
$$
Thus
$$
\left\|T_2 f\right\|_{L_{\varphi,\,w}}\le C\,\|f\|_{L_{p,v}(R^n)}\eqno(4.13)
$$
where $C> 0$ is independent of $f$ and $x\in R^n.$

Combining the inequalities (4.11),(4.12) and (4.13) we obtain the proof of Theorem 5.

In particular, for $\varphi(x, t)= t^{q(x)}$ by virtue of Remark 2 we have the following Corollary.

\begin{corollary}{[2]} Let $1< p\le q(x)\le \overline q< \infty$ and $x\in R^n.$ Suppose that
$v(x)$ and $w(x)$ are weight functions defined on $R^n$ and satisfying conditions (4.8), (4.9) and (4.10) of Theorem 5. Let $T$ be a sublinear operator satisfies the condition (4.7) of Theorem 5.

Then there exists a positive constant $C,$ independent of $f$, such
that for all $f\in L_{p, v}(R^n)$
$$
\|Tf\|_{L_{q(\cdot), w}\left(R^n\right)}\le C \|f\|_{L_{p, v}\left(R^n\right)}.
$$
\end{corollary}

\begin{remark} Note that the condition (4.7) was introduced in [44].
Many interesting operators in harmonic analysis,
such as the Calderon-Zigmund singular integral operators,
Hardy-Littlewood maximal operators, Fefferman's singular integrals,
Ricci-Stein's oscillatory singu-lar integrals, Bochner-Riesz means
and so on is satisfied the condition (4.7). In the case $p(x)= p=
const$ for classical Lebesgue spaces the Theorem 5 was proved in
[47] (see also [18] and [33]). Also, for classical Lebesgue spaces in
[16] and [19] was found new type sufficient conditions on weights for Calderon-Zigmund singular
integral operator, whenever the weight functions are radial monotone
functions. In particular, the boundedness of certain convolution operator
in a weighted Lebesgue space with kernel satisfying the generalized H\"{o}rmander's
condition was proved in [5].
\end{remark}

{\bf Acknowledgement.} Author was supported by the Science Development Foundation under the President of
the Republic of Azerbaijan EIF-2010-1(1)-40/06-1.

\vspace{6mm}

\begin{center}
{\bf References}
\end{center}

[1] R.A.Bandaliev, {\it On an inequality in Lebesgue space with mixed norm and with variable summability
exponent},  Matem. Zametki, {\bf 84}(2008), no. 3, 323-333 (in Russian): English translation: in
Math. Notes, {\bf 84}(2008), no. 3, 303-313. MR2473750 (2010k:46030).

[2] R.A.Bandaliev,{\it The boundedness of certain sublinear operator in the weighted variable Lebesgue spaces},
Czechoslovak Math. J. {\bf 60}(2010), no. 2, 327-337. MR2657952 (2011g:47068)

[3] R.A.Bandaliev, {\it The boundedness of multidimensional Hardy operator in the weighted variable
Lebesgue spaces}, Lithuanian Math. J. {\bf 50}(2010), no.3, 249-259. MR2719561 (2011h: 42014)

[4] R.A.Bandaliev, {\it Criteria of two-weighted inequalities for multidimensional
Hardy type operator in weighted Musielak-Orlicz spaces and some application},  Math. Nachr. 2012 (accepted).

[5] R.A.Bandaliev and K. K. Omarova, {\it Two-weight norm inequalities for certain
singular integrals}, Taiwanese Jour. of Math. {\bf 16}(2012), no. 2, 713-732.

[6] C. Bennett and R. Sharpley, {\it Interpolation of operators}, Pure
Appl. Math. 129, Academic Press, 1988. MR 0928802 (89e:46001)

[7] E. I. Berezhnoi, {\it Sharp estimates of operators on the cones of ideal
spaces}, Trudy Mat. Inst. Steklov. {\bf 204}(1993), 3-36; English transl. in Proc. Steklov Inst. Math. {\bf 204}(1994), no. 3 . MR 96a:46052

[8] E. I. Berezhnoi,{\it Two-weighted estimations for the Hardy-Littlewood maximal function in ideal Banach
spaces}, Proc. Amer. Math. Soc. {\bf 127}(1999), 79-87. MR1622773 (99d:42028)

[9] S.Bloom and R.Kerman,{\it Weighted $L_\Phi$ integral inequalities for operators of Hardy type},
Studia Math. {\bf 110}(1994), no. 1,  35-52. MR 95f:42031

[10] J.Bradley,{\it Hardy inequalities with mixed norms}, Canadian Math. Bull. {\bf 21}(1978), 405-408. MR 80a:26005

[11] A. Cianchi, {\it Hardy inequalities in Orlicz spaces}, Trans. Amer. Math. Soc.{\bf 351}(1999), 2459-2478. MR1433113 (99i:46016)

[12] J.A.Cochran and C.-S.Lee, {\it Inequalities related to Hardy's and Heinig's}, Math. Proc. Camb. Phil.Soc.
{\bf 96}(1984), 1-7. MR 86g:26026

[13] D.Cruz-Uribe, A.Fiorenza and C.J.Neugebauer, {\it The maximal
function on variable $L^p$ spaces},  Ann. Acad. Sci. Fenn. Math. {\bf
28}(2003), no. 1, 223-238; corrections in Ann. Acad. Sci.
Fenn. Math. {\bf 29}(2004), no. 2, 247-249. MR 1976842 (2004c:42039) $\,$ MR 2041952 (2004m:42018)

[14] L.Diening and S.Samko, {\it Hardy inequality in variable exponent
Lebesgue spaces}, Frac. Calc. Appl. Anal.{\bf 10}(2007), no. 1, 1-18.

[15] L.Diening, P.Harjulehto, P.H\"{a}st\"{o} and M. R\.{u}\v{z}i\v{c}ka,
{\it Lebesgue and Sobolev spaces with variable exponents,} Springer Lecture Notes,
v.2017, Springer-Verlag, Berlin, 2011.

[16] D.E.Edmunds and V.Kokilashvili, {\it Two-weighted inequalities
for singular integrals}, Canadian Math.Bull. {\bf 38}(1995), 295-303.

[17] D.E.Edmunds, V.Kokilashvili and A.Meskhi,{\it On the boundedness and
compactness of weighted Hardy operators in spaces $L^{p(x)}$}, Georgian Math.J. {\bf 12}(2005),
no. 1, 27-44. MR 2136721

[18] A.D.Gadjiev and I.A.Aliev,{\it Weighted estimates for
multidimensional singular integrals generated by a generalized shift
operator}, Mat. Sbornik {\bf 183}(1992), no. 9, 45-66.(In Russian.). English transl. in Sb. Math.
{\bf 77}(1994), no. 1, 37-55. MR1198834(94h:41046)

[19] V.S.Guliev, {\it Two-weight inequalities for integral operators in $L_p$-spaces and their applications},
Trudy Mat. Inst. Steklov. {\bf 204}(1993), 113-116; English transl. in Proc. Steklov Inst. Math. {\bf 204}(1994), 97-116.  MR1320021 (96j:26029)

[20] G.H.Hardy, J.E. Littlewood, and G.P\'{o}lya, {\it Inequalities}, Cambridge Univ. Press, 1988. MR944909 (89d:26016)

[21] P.Harjulehto, P.H\"{a}st\"{o} and M.Koskenoja, {\it Hardy's inequality in a variable exponent
Sobolev space},  Georgian Math.J. {\bf 12}(2005), no. 1, 431-442.

[22] H.P.Heinig, {\it Some extensions of  inequalities}, SIAM J. Math. Anal. {\bf 6}(1975), 698-713.

[23] P.Jain, L.E.Persson and A.Wedestig, {\it From Hardy to Carleman and general mean-type inequalities},
Function Spaces and Appl., CRC Press (New York)/Narosa Publishing House
(New Delhi)/Alpha Science (Pangbourne), (2000), 117-130.

[24] P.Jain, L.E.Persson and A.Wedestig, {\it Carleman-Knopp type inequalities via Hardy's inequality},
Math. Inequal. Appl. {\bf 4}(2001), no. 3, 343-355.

[25] K.Knopp,  {\it \"{U}ber Reihen mit positiven Gliedern}, J. London Math. Soc. {\bf 3}(1928), 205-211.

[26]  T.S.Kopaliani, {\it On some structural properties of Banach function
spaces and bounded-ness of certain integral operators}, Czechoslovak Math. J.
({\bf 54}){\bf 129}(2004), 791-805.

[27] M.Krbec, B.Opic, L.Pick and J.R\'{a}kosn\'{\i}k, {\it Some recent
results on Hardy type operators in weighted function spaces and
related topics}, in: "Function Spaces, Differential Operators and
Nonlinear Analysis", Teubner, Stuttgart,(1993), 158-184.

[28] A.Kufner and L.E.Persson, {\it Integral inequalities with weights} (World Scientific Publish-ers, Singapore), 2002.

[29] J.Lindenstrauss and L.Tzafriri, {\it Classical Banach spaces II}, Ergebnisse der Mathematik und ihrer
Grenzgebiete 97, Springer-Verlag, Berlin-Heidelberg-New York, 1979.  MR 81c:46001

[30] E.Lomakina and V.Stepanov, {\it On the Hardy-type integral operator in Banach function spaces}, Publicacions Matem\`{a}tiques, {\bf 42}(1998), 165-194 .

[31] E.R. Love, {\it Inequalities related to those of Hardy and of Cochran and Lee},
Math. Proc. Camb. Phil.Soc. {\bf 99}(1986), 395-408.

[32] W.A.J.Luxemburg,{\it Banach function spaces}, Thesis, Delfi, 1955.

[33] G.Lu, Sh.Lu and D.Yang, {\it Singular integrals and
commutators on homogeneous groups}, Analysis Math. {\it 28}(2002),
103-134.

[34] F.I.Mamedov and A.Harman,  On a weighted inequality of Hardy
type in spaces $L^{p(\cdot)},$ J. Math. Anal. Appl. {\bf 353}(2009), no. 2, 521-530.

[35] R.A.Mashiyev, B.\c{C}eki\c{c}, F.I.Mamedov and S.Ogras,{\it Hardy's inequality in power-type weighted
$L^{p(\cdot)}(0, \infty)$ spaces}, J. Math. Anal. Appl. (1) {\bf 334}(2007), no. 1, 289-298.

[36] V.G.Maz'ya,  {\it Sobolev spaces}, (Springer-Verlag, Berlin, (1985). MR 87g:46056

[37] B.Muckenhoupt, {\it Hardy's inequality with weights},  Studia Math. {\bf 44}(1972), 31-38. MR 47:418

[38] J.Musielak, {\it Orlicz spaces and modular spaces},  Lecture
Notes in Math. (1034), Springer-Verlag, Berlin-Heidelberg-New York, 1983. MR 724434 (85m:46028)

[39] B.Opic and P.Gurka, {\it Weighted inequalities for geometric means},
Proc. Amer. Math. Soc. {\bf 120}(1994), no. 3, 771-779. MR1169043(94e: 26036)

[40] L.-E. Persson and V.D. Stepanov, {\it Weighted integral inequalities with the
geometric mean operator}, J. Inequal. Appl. {\bf 7}(2002), no. 5, 727-746.

[41] L.Pick and B.Opic,  {\it On geometric mean operator}, J. Math. Anal. Appl. {\bf 183}(1994), no. 3,
652-662.

[42] L.Quinsheng, {\it Two weight $\Phi$-inequalities for the Hardy operator, Hardy-Littlewood maximal
operator and fractional integrals}, Proc. Amer. Math. Soc. {\bf 118}(1993), no. 1, 129-142. MR1279455(99j: 36037)

[43] A.Schep, {\it Minkowski's integral inequality for function norms}, Operator theory. Operator theory
in function spaces and Banach lattices. {\bf 75}(1995), Oper. Theory Adv. Appl., 299-308.

[44] F.Soria and G.Weiss, {\it A remark on singular integrals and
power weights}, Indiana Univ.Math.J. {\it 43}(1994), 187-204. MR1275458

[45] G.Tomaselli, {\it A class of inequalities},  Boll.Un.Mat.Ital. {\bf 2}(1969), 622-631.

[46] A. Wedestig, {\it Some new Hardy type inequalities and their limiting inequalities}, J.of Ineq. in Pure
and  Appl. Math. {\bf 61}(2003), no. 4, 1-33.

[47] Y.Zeren and V.S.Guliyev, {\it Two-weight norm inequalities for some anisotropic sublinear operators},
Turkish Math.J. {\bf 30}(2006), 329-350.

\end{document}